\newcommand{\CA}{\hbox{{$\cal A$}}}

\newcommand{\Z}{\mathbb{Z}}
\newcommand{\C}{\mathbb{C}}

\newcommand{\note}[1]{}
\newcommand{\h}{{\scriptstyle\frac{1}{2}}}
\newcommand{\extd}{{\rm d}}
\newcommand{\image}{{\rm image}}
\newcommand{\isom}{{\cong}}
\newcommand{\eps}{{\epsilon}}
\newcommand{\tens}{\mathop{\otimes}}
\newcommand{\id}{{\rm id}}

\newcommand{\la}{\triangleright}
\renewcommand{\o}{{}_{\scriptscriptstyle(1)}}
\renewcommand{\t}{{}_{\scriptscriptstyle(2)}}
\newcommand{\bo}{{}^{\bar{\scriptscriptstyle(1)}}}
\newcommand{\bt}{{}^{\bar{\scriptscriptstyle(2)}}}

\newcommand{\proof}{{\bf Proof\ }}
\newcommand{\eproof}{$\quad \diamond$\bigskip}

\newcommand{\und}[1]{{\underline {#1}}}

\documentclass[11pt,A4]{article}
\usepackage{amssymb,amsmath,epsfig}
\newtheorem{lemma}{Lemma}[section]
\newtheorem{propos}[lemma]{Proposition}

\newtheorem{corol}[lemma]{Corollary}

\begin{document}\baselineskip 20pt

{\ }\qquad \hskip 4.3in \vspace{.2in}

\begin{center} {\LARGE NONCOMMUTATIVE COHOMOLOGY AND ELECTROMAGNETISM
ON $\C_q[SL_2]$ AT ROOTS OF UNITY}
\\ \baselineskip 13pt{\ }\\
{\ }\\   Xavier Gomez\footnote{EU Human Potential postdoc} + Shahn
Majid\footnote{Professor and Royal Society University Research
Fellow} \\ {\ }\\ School of
Mathematical Sciences\\
Queen Mary, University of London, Mile End Rd\\ London E1 4NS, UK
\end{center}
\begin{center}
October 2001
\end{center}

\begin{quote}\baselineskip 14pt
\noindent{\bf Abstract} We compute the noncommutative de Rham
cohomology for the finite-dimensional q-deformed coordinate ring
$\C_q[SL_2]$ at odd roots of unity and with its standard
4-dimensional differential structure. We find that $H^1$ and $H^3$
have three additional modes beyond the generic $q$-case where they
are 1-dimensional, while $H^2$ has six additional modes. We solve
the spin-0 and Maxwell theory on $\C_q[SL_2]$ including a complete
picture of the self-dual and anti-self dual solutions and of
Lorentz and temporal gauge fixing. The system behaves in fact like
a noncompact space with self-propagating modes (i.e., in the
absence of sources). We also solve with examples of `electric' and
`magnetic' sources including the biinvariant element $\theta\in
H^1$ which we find can be viewed as a source in the local
(Minkowski) time-direction (i.e. a uniform electric charge
density).

\bigskip
{\em Keywords:} noncommutative geometry, roots of unity, quantum
groups, cohomology, electromagnetism, light
\end{quote}

\baselineskip 20pt
\section{Introduction}

By now there is a standard formulation of differential calculi or
`exterior algebra' of differential forms on quantum groups such as
$\C_q[SL_2]$. The standard bicovariant ones correspond essentially
to representations\cite{Ma:cla}, i.e. are labelled in this case by
spin $j\in\h \Z_+$ and have dimension $(2j+1)^2$ (there are also
exotic twists of the standard ones which do not concern us). In
our case the smallest nontrivial calculus is 4 dimensional and was
already known since the earliest works \cite{Wor:dif}. The entire
exterior algebra and exterior derivative are also known, and it is
known that dimensions in each degree of forms and the resulting
cohomology for generic $q$ are\cite{Gri:bic}
\[\dim(\Omega)=1:4:6:4:1,\quad H^0=\C,\quad H^1=\C,\quad H^2=0,
\quad H^3=\C,\quad H^4=\C.\] The nontrivial generator in degree 1
is the bi-invariant element $\theta$ that defines $\extd$ by
graded-commutator. The further physics and geometry on such spaces
has been mainly looked at for generic $q$, where (with some
modifications such as a 1-dimensional extension) it follows
broadly the line of the undeformed case.

What we show in the present purely computational paper is the
existence of completely different and novel phenomenona when,
however, $q$ is an odd root of unity. This case is in many ways
more relevant to both physics (e.g in the Wess-Zumino-Witten
model) and mathematics (e.g. the image of the quantum Frobenius
map and because of known links to group theory in finite
characteristic). We work with the reduced finite-dimensional
quantum group, which is then a nonsemisimple Hopf algebra. This
has the merit that all linear (and some nonlinear) aspects of the
geometry can be fully computed. The model also contrasts markedly
from the case of finite group algebras recently studied
elsewhere\cite{MaRai:ele}. For the differential calculus itself
the theorem for factorisable quantum groups in \cite{Ma:cla}
implies that these are classified by two-sided ideals in
$u_q(sl_2)$. So the smallest nontrivial calculus is again the
4-dimensional one, which is the calculus that we use. Its
structure is recalled briefly in Section~2.

We then find in Section~3 that there are additional elements of
$H^i$ not present for generic $q$. In all cases that we have
checked (namely 3,5,7'th roots) we find in fact that
\[ H^\cdot\isom \Lambda,\]
the space of right-invariant forms as a graded vector space, as
well as an exact sequence for the $H^i$. The additional cohomology
modes correspond to topological gauge fields with zero curvature
in the Maxwell theory reflecting nontrivial topology created by
the quotienting to the reduced quantum group. In Section~4 we use
the family of
 `Killing form' metrics in \cite{Ma:rieq} and show how the requirement
 of $\star^2=\id$ for the Hodge-*
operator singles out a particular $q$-deformed Minkowski one (this
applies for generic $q$). We then proceed to solve the Maxwell
theory for $r=3$ completely. Among interesting features, we find
that for spin 0 the wave operator $\square$ is not fully
diagonalisable (this is due to the nonsemisimplicity), while on
the other hand every solution of the sourceless Maxwell's
equations may be written as a sum of a self-dual and an
antiself-dual solution. We are also able to completely analyse
gauge fixing issues which are usually glossed over in gauge
theories in physics; we find the novel result that not all
solutions can be rendered in Lorentz gauge, nor all in temporal
gauge, but that the two gauges between them `patch' the moduli of
solutions. We expect the phenomena found here by computation to
hold for all odd roots.

\section{Exterior algebra}

Here we fix the algebras and exterior algebras in question in
notation that we will use. In effect, in order to have reliable
formulae for root of unity we carefully compute the (well-known)
4-D calculus from a modern crossed-module point of view. We let
$q^2\ne 1$. The quantum group $\CA=\C_q[SL_2]$ has a matrix of
generators $t^i{}_j=\begin{pmatrix}a& b\\ c& d\end{pmatrix}$ with
relations
\[ ba=qab,\quad ca=qac,\quad db=qbd,\quad dc=qcd,\quad cb=bc,
\quad da-ad=q\mu bc,\quad ad-q^{-1}bc=1,\] where $\mu=1-q^{-2}$,
and the matrix coalgebra structure. For its 4d calculus we take a
basis $e_i{}^j=\begin{pmatrix}e_a& e_b\\ e_c & e_d\end{pmatrix}$,
where $e_1{}^2=e_b$, of the space $\Lambda^1$ of right-invariant
differential 1-forms. This space $\Lambda^1$ is specified as a
left $\C_q[SL_2]$-crossed module, namely
 with coaction and action
\[ \Delta_L e_i{}^j=(S t^k{}_i)t^j{}_l\tens e_k{}^l,\quad a\la
\begin{pmatrix}e_a&e_b\\
e_c&e_d\end{pmatrix}=\begin{pmatrix} qe_a+q\mu^2 e_d& e_b\\
e_c& q^{-1}e_d\end{pmatrix}\]
\[ b\la \begin{pmatrix}e_a&e_b\\ e_c&e_d\end{pmatrix}
=\begin{pmatrix}\mu e_c & q\mu e_d\\ 0&0 \end{pmatrix},\quad
 c\la \begin{pmatrix}e_a&e_b\\ e_c&e_d\end{pmatrix}
 =\begin{pmatrix}\mu e_b& 0\\ q\mu e_d&0 \end{pmatrix},\quad
 d\la \begin{pmatrix}e_a&e_b\\ e_c&e_d\end{pmatrix}
 =\begin{pmatrix}q^{-1}e_a& e_b\\ e_c& q e_d \end{pmatrix}.\]
$\Omega^1=\Lambda^1\tens \CA$ is generated by these forms as a
free right module over $\C_q[SL_2]$ while as a bimodule the left
action is $h.e_a=(h\o\la e_a) h\t$ for all $h\in\C_q[SL_2]$, etc.
This comes out as
\begin{eqnarray*}
&&[\begin{pmatrix}c\\ d\end{pmatrix},e_b]=[\begin{pmatrix}a\\
b\end{pmatrix},e_c]=[\begin{pmatrix}a\\ b
\end{pmatrix},e_d]_{q^{-1}}=[\begin{pmatrix}c\\
d\end{pmatrix},e_d]_q=0,\\
&&[\begin{pmatrix}a\\ b\end{pmatrix},e_b]=q\mu e_d
\begin{pmatrix}c\\ d\end{pmatrix},\qquad [\begin{pmatrix}c\\ d
\end{pmatrix},e_c]=q\mu e_d \begin{pmatrix}a\\ b\end{pmatrix},\\
&&[\begin{pmatrix}c\\
d\end{pmatrix},e_a]_{q^{-1}}=\mu e_b \begin{pmatrix}a\\
b\end{pmatrix}, \quad  [\begin{pmatrix}a\\
b\end{pmatrix},e_a]_q=\mu e_c \begin{pmatrix}c\\ d\end{pmatrix}
+q\mu^2 e_d \begin{pmatrix}a\\ b\end{pmatrix},\end{eqnarray*}where
$[x,y]_q=xy-qyx$.

Also from the crossed module structure is the braiding
$\Psi(e_a\tens e_b)=e_a\bo\la e_b\tens e_a\bt$, etc., where
${}\bo$ and ${}\bt$ denote the outputs of $\Delta_L$. This comes
out as
\begin{eqnarray*}
\Psi(e_a \otimes e_a ) &=& e_a\otimes e_a -\mu (e_b\otimes
e_c-e_c\otimes e_b)+
    q\,\mu^2\, e_d\otimes (q\, e_a-q^{-1} e_d)
\\
\Psi(e_b \otimes e_b ) &=& e_b\otimes e_b
\\
\Psi(e_c \otimes e_c ) &=& e_c\otimes e_c
\\
\Psi(e_d \otimes e_d ) &=& e_d\otimes e_d
 \\
\Psi(e_a \otimes e_d ) &=& e_d\otimes e_a
\\
\end{eqnarray*}
\begin{eqnarray*}
\Psi(e_d \otimes e_a ) &=& e_a\otimes e_d+\mu\, (e_b\otimes
e_c-e_c\otimes e_b)
    - q\,\mu^2\, e_d\otimes (q\, e_a-q^{-1} e_d)
\\
\Psi(e_b \otimes e_c ) &=& e_c\otimes e_b+ q\, \mu\, e_d\otimes
(q\, e_a-q^{-1} e_d)
\\
\Psi(e_c \otimes e_b ) &=& e_b\otimes e_c-q\, \mu\, e_d\otimes
(q\, e_a-q^{-1} e_d)
\\
\Psi(e_a \otimes e_b ) &=& e_b\otimes e_a+q^2\, \mu\, e_d\otimes
e_b
\\
\Psi(e_b \otimes e_a) &=& q^{-2}\, e_a\otimes e_b+
    \mu\,q^{-1} e_b\otimes (q\, e_a- q^{-1}e_d)\\
\Psi(e_a\otimes e_c) &=& e_c\otimes e_a-\mu\, e_d\otimes e_c
\\
\Psi(e_c\otimes e_a) &=& q^2\, e_a\otimes e_c-
    \mu\,q\, e_c\otimes (q\, e_a- q^{-1}e_d)
    +[2]_{q^2}\, \mu^2\, e_d\otimes e_c
\\
\Psi(e_b \otimes e_d ) &=& q^2\, e_d\otimes e_b
\\
\Psi(e_d \otimes e_b ) &=& e_b\otimes e_d- q^2\, \mu\, e_d\otimes
e_b
\\
\Psi(e_c \otimes e_d ) &=& q^{-2}\, e_d\otimes e_c
\\
\Psi(e_d \otimes e_c ) &=& e_c\otimes e_d+\mu\, e_d\otimes e_c
\end{eqnarray*}
where $[n]_q=(1-q^n)/(1-q)$.

This extends as a bimodule map to an endomorphism of
$\Omega^1\tens_{\C_q[SL_2]}\Omega^1$. Following Woronowicz we then
define $\Omega^2=\Omega^1\tens_M\Omega^1/\ker(\id-\Psi)$, etc.
Equivalently, in a modern braided group
approach\cite{Ma:eps}\cite{BesDra:dif} which is computationally
easier, $\Omega$ is a free right $\C_q[SL_2]$ module over the
invariant exterior forms
\[ \Lambda=T\Lambda^1/\oplus_n\ker A_n;\quad A_n=[n,-\Psi]!
=(\id\tens A_{n-1})[n,-\Psi] \]
\[[n,-\Psi]=\id-\Psi_{12}+\Psi_{12}\Psi_{23}\cdots
+(-1)^{n-1}\Psi_{12}\cdots \Psi_{n-1,n}.\] Here $[n,-\Psi]$ are
the braided-integers\cite{Ma:book} induced by a braiding $-\Psi$
and $\Lambda$ is a braided group with additive coproduct
$\und\Delta e_a=e_a\tens 1 +1\tens e_a$, etc. The above relations
ensure that it is dually paired with a similar braided group
$\Lambda^*$ and these together ensure Poincar\'e duality. In
particular, $A_2=\id-\Psi$ and hence the relations in degree 2,
which are in fact all the relations for generic $q$, come out as:
\[ <e_b,e_c,e_d>\quad {\rm usual\ Grassmann\ algebra},\quad e_a^2
=\mu e_b\wedge e_c,\quad e_a\wedge e_d+e_d\wedge e_a+\mu e_b\wedge
e_c=0,\]
\[e_a\wedge e_b+q^2e_b\wedge e_a-\mu e_b\wedge e_d=0,\quad
e_c\wedge e_a+q^2e_a\wedge e_c+\mu e_c\wedge e_d=0\]

Note that if we define the corresponding symmetric algebra by
$T\Lambda^1/\image(\id+\Psi)$ then we have q-Minkowski space in
the braided-matrix form. The exterior algebra in this case has a
similar form\cite{Ma:book} to the above in terms of exact
differentials, since both come from Meyer's braiding for the
additive braided group structure of q-Minkowski space.

Finally, the exterior derivative is
\[ \extd =-[\theta,\ \},\quad \theta=e_a+e_d\]
where we use the commutator on even degree $n$ and anticommutator
on odd. Note that in our conventions $\bar\extd=[\ ,\theta \}$ is
more natural but would be a right-derivation. The element $\theta$
is closed but not exact and is biinvariant. Explicitly,
\[ \extd \begin{pmatrix}a\\ b\end{pmatrix}=(q-1)(e_a-q^{-1}
(1-\mu[2]_q) e_d)\begin{pmatrix}a\\ b\end{pmatrix} +\mu
e_c\begin{pmatrix}c\\ d\end{pmatrix},\quad  \extd
\begin{pmatrix}c\\ d\end{pmatrix}=(q-1)(e_d-q^{-1}e_a)
\begin{pmatrix}c\\ d\end{pmatrix}+\mu e_b\begin{pmatrix}a\\ b\end{pmatrix}\]
\[ \extd e_a=-\mu e_b\wedge e_c,\quad \extd e_d=\mu e_b\wedge e_c,\quad
\extd e_b=-\mu(e_a\wedge e_b+q^{-2}e_b\wedge e_d),\quad \extd e_c
=\mu(q^{2}e_a\wedge e_c+e_c\wedge e_d).\]

\section{Roots of unity and cohomology}

We now study $\CA=\C_q[SL_2]$ reduced at $q^r=1$ a primitive
$r$'th root of unity by the additional relations
\[ c^r=b^r=0,\quad a^r=d^r=1,\]
which we suppose from now on. Here $a^r,b^r,c^r,d^r$ generate an
undeformed $\C[SL_2]$ central sub-Hopf algebra of the original
$\C_q[SL_2]$. Note also that in the reduced case
$d=a^{-1}(1+q^{-1}bc)$ is redundant and moreover the algebra
becomes finite dimensional, with $\dim(\CA)=r^3$. A basis of $\CA$
is $\{a^mb^nc^k\}$ for $0\le m,n,k\le r-1$. All kernel
computations are done below for $r=3,5,7$ for concreteness, but we
expect identical results for all odd $r$.

\begin{propos} At least for $r=3,5,7$, the exterior algebra for the
reduced quantum group has the same dimensions as for generic $q$
(namely 1:4:6:4:1) and is given entirely by relations in degree 2
(a quadratic algebra). Moreover, the exterior derivative descends
to one over the reduced $\C_q[SL_2]$.
\end{propos}
\proof Since the reduced quantum group remains a Hopf algebra is
bicovariant calculi are still defined by quotient crossed modules
of $\ker\eps$. Our particular crossed module remains one with the
same form of action and coaction, bimodule structure and braiding.
Therefore it is only a matter of computing the explicit braided
factorial matrices $[n,-\Psi]!$ for $n=2,3,4$ and in particular
the dimensions of their kernel, which we find to be the same
provided $r$ is odd (for example $r=6$ is different). Hence the
algebra $\Lambda$ is unchanged in this case. Since we have not
discussed explicity the projection from $\ker\eps$ to $\Omega_0$
we also verify directly that $\extd$ is consistent with the
additional relations of the reduced quantum group. \eproof.

Next, in order to compute cohomology we need $\extd$ on a general
element of $\CA$. This is given by the Leibniz rule and the
following:
\begin{lemma} For all invertible $q^2\ne 1$,
\begin{eqnarray*}
\extd(a^m b^n c^k)&=&  \quad e_a\, . (q^{m+n-k}-1)\, a^m b^n c^k
+\mu\, e_b\, .\, q^{n-k+1}[k]_{q^2}\, a^{m+1}\, b^n\, c^{k-1}
\\
    && +\mu\, e_c\, .\, q^{-k-n}\,
(\, [m+n]_{q^2}\, a^{m-1}b^nc^{k+1}+q[n]_{q^2}\,
a^{m-1}b^{n-1}c^k)
\\
&& +\mu^2\, e_d\, .\, q^{-k-m-n+2}(\, [k+1]_{q^2}\, [m+n]_{q^2}\,
a^m b^n c^k+
        q[n]_{q^2}\, [k]_{q^2}\, a^mb^{n-1}c^{k-1})
\\
    && +e_d\, .\, (q^{-m-n+k}-1)\, a^m b^n c^k
\end{eqnarray*}
 \end{lemma}
\proof  We first iterate the stated bimodule relations to obtain
\begin{eqnarray*} u^n.\,e_a &=& q^{+n}\, e_a.\,u^n+q^{n-1}\, [n]\,
\mu\, e_c.\, u^{n-1}v
    +q\mu^2\, [n]\, e_d.\,u^n
\\
v^n.\,e_a &=& q^{-n}\, e_a.\,v^n+q^{1-n}\, [n]\, \mu\, e_b.\,
v^{n-1}u
\\ u^n.\,e_b &=& e_b.\,u^n+q\mu\, [n]\, e_d.\, u^{n-1}v
\\
v^n.\,e_b &=& e_b.\,v^n
\\ u^n.\,e_c &=& e_c.\,u^n
\\
v^n.\,e_c &=& e_c.\,v^n+q\mu\, [n]\, e_d.\, v^{n-1}u
\\ u^n.\,e_d &=& q^{-n}e_d.\,u^n
\\
v^n.\, e_d &=& q^{+n}\, e_d.\,v^n
\end{eqnarray*}
where $[n]=(q^n-q^{-n})/(q-q^{-1})$ and $u=\begin{pmatrix}a\\
b\end{pmatrix}$, $v=\begin{pmatrix}c\\ d\end{pmatrix}$. Then by
recurrence, one gets, for $X=a^mb^nc^k$,
\begin{eqnarray*}
X\, .\, e_a &=& \quad e_a\, . q^{m+n-k}\,\, X +\mu\, e_b\, .\,
q^n[k]\, a^{m+1}\, b^n\, c^{k+1}
\\
    && +\mu\, e_c\, .\, q^{n-1}\,
        (\, q^{m-k}[m]\, a^{m-1}b^nc^{k+1}+[n]\, a^mb^{n-1}c^kd)
\\
    && +q\mu^2\, e_d\, .\, (\, [k+1]\, [m+n]\, X+
        q^{-m}[n]\, [k]\, a^mb^{n-1}c^{k-1})
\\ X\, .\, e_d &=& e_d\, .\, q^{-m-n+k}\, X
\end{eqnarray*}
We then compute $\extd X=X\theta-\theta X$. \eproof

Finally, we choose an explicit basis for each degree of the
exterior algebra. Here $e_{abc}\equiv e_a\wedge e_b\wedge e_c$,
etc. for our chosen basis elements
\[ \Lambda^2=\{e_{ab},e_{ac},e_{ad},e_{bc},e_{bd},e_{cd}\},
\quad \Lambda^3=\{e_{abc},e_{abd},e_{acd},e_{bcd}\},
\quad \Lambda^4=\{e_{abcd}\}\] and we then use the above relations
to explicitly define $\wedge,\extd$ on right-invariant forms as a
$16\times 6$ matrix and a $4\times 6$ matrix respectively. With
these ingredients it is a matter of linear algebra to compute
cohomology.

\begin{propos} At least for $r=3,5,7$ the noncommutative de
Rham cohomology $H^1$ for the $4d$ calculus on the reduced quantum
group $\C_q[SL_2]$ is 4 dimensional with basis
\[ \theta=e_a+e_d,\quad h_1=e_b ac^{r-1},\quad h_2=e_ca^{r-1}b^{r-1},
\quad n=e_a+e_ca^{r-1}c\]
\end{propos}
\proof We write $\extd_0:\Omega^0\to\Omega^1$ as an $r^3\times
4r^3$ matrix. We also compute the wedge product
$\Lambda^1\tens\Lambda^1\to \Lambda^2$ and exterior derivative
$\extd:\Lambda^1\to\Lambda^2$ as explained above. The Leibniz rule
then allows us to define $\extd_1:\Omega^1\to\Omega^2$ from these
ingredients as a $4r^3\times 6r^3$ matrix. We then compute the
null spaces to find the dimension of the cohomology (for example
the kernel of $\extd_1$ is 30-dimensional for $r=3$ and
346-dimensional for $r=7$, while the image of $\extd_0$
26-dimensional and 342-dimensional respectively). We also verify
$\extd_1\extd_0=0$ as a programming check. Finally we chose 4
vectors in $\ker\extd_1$ and not in the image of $\extd_0$ and
verify that together with a basis of the image of $\extd_0$ they
form a linearly independent set, i.e. their classes provide a
basis of $H^1$. \eproof

We find the same kind of phenomenon for the higher cohomology.

\begin{propos} At least for $r=3,5,7$ the  cohomologies have the
same dimensions as $\Lambda$ in each degree. As representatives
we have:
\begin{eqnarray*} H^2:&&  m_1=e_{bd}ac^{r-1},
\quad m_2=e_{ab}ac^{r-1},
\quad m_3=e_{ac}a^{r-1}b^{r-1}, \quad  m_4=e_{cd}a^{r-1}b^{r-1},\\
&& m_5=(e_{ac}-e_{cd})a^{r-1}c-e_{ad},\quad
m_6=e_{bd}ab^{r-1}c^{r-2}+q^4 e_{cd}a^{r-1}b^{r-2}c^{r-1}.\\
 H^3:&& \Theta=e_{bcd}b^{r-1}c^{r-1},\quad
 h_1^*=e_{abd}ac^{r-1},\quad
 h_2^*=e_{acd}a^{r-1}b^{r-1},
\\
&& s=e_{abd}ab^{r-1}c^{r-2}+q^4
e_{acd}a^{r-1}b^{r-2}c^{r-1},\\
H^4:&& e_{abcd}b^{r-1}c^{r-1}.\end{eqnarray*}
 \end{propos} \proof
We proceed with respect to our basis above to compute the wedge
products $\Lambda^2\tens\Lambda^1\to \Lambda^3$ and
$\Lambda^1\tens\Lambda^2\to\Lambda^3$ as $24\times 4$ matrices. We
also use the (graded) Leibniz rule to define
$\extd:\Lambda^2\to\Lambda^3$ using these projectors and the
matrices already computed for $H^1$. Finally we combine these via
Leibniz with $\extd_0$ to obtain $\extd_2:\Omega^2\to \Omega^3$ as
a $6r^3\times 4r^3$ matrix, and compute its kernel and the image
of $\extd_1$ above (for example the kernel of $\extd_2$ is
84-dimensional for $r=3$ and 1032-dimensional for $r=7$, while the
image of $\extd_1$ is $78$-dimensional and 1026-dimensional
respectively). Similarly we proceed to
$\extd_3:\Omega^3\to\Omega^4$ (its kernel is 82-dimensional for
$r=3$ and $1030$-dimensional for $r=7$, while the image of
$\extd_2$ is 78-dimensional and 1026-dimensional respectively).
The image of $\extd_3$ has codimension 1, so $H^4$ is similarly
1-dimensional. We then chose representatives and check linear
independence in the quotient spaces. Our notations for them will
be relevant later. Note also that the kernel of $\extd_0$ is
1-dimensional and $H^0$ clearly has a basis given by 1. \eproof

Finally, we observe that the cohomology is itself a complex under
the operation $\theta\wedge$ since $\theta\wedge\theta=0$ and
$\extd\theta=0$.

\begin{propos} At least for $r=3,5,7$, the sequence
$0\to H^0\to H^1\to H^2\to H^3\to H^4\to 0$
defined by $\theta\wedge$ is exact.
\end{propos}
\proof We let $\theta_0$ be $\theta$ acting by multiplication in
degree 0, etc. The image of $\theta_0$ is $\theta$. Its complement
in the basis shown has image
\[\theta\wedge h_1=m_2-m_1,\quad \theta \wedge h_2
=m_3-m_4,\quad \theta\wedge n=m_5\] which is a 3-dimensional
subspace of $H^2$. Its complement has basis $m_1+m_2$, $m_3+m_4$,
$m_6$. These map under $\theta_2$ up to normalisation to
$h_1^*,h_2^*,s$, which are three of the basis elements of $H^3$.
Their complement $\Theta$ maps under $\theta_3$ to the generator
of $H^4$. \eproof

\section{Wave equations and Hodge-*}

Next we describe the Hodge * operator corresponding to the
`Killing metric' introduced in \cite{Ma:rieq}.  These are further
geometric structures on the full $\C_q[SL_2]$ and after recalling
them in the form that we need, we will then specialise to our
reduced root of unity case. In our conventions the general metric
is:

\begin{lemma} For all invertible $q^2\ne 1$,
\[\eta\equiv \eta^{ij}e_i\tens e_j=e_b\tens e_c
+q^2 e_c\tens e_b+{(qe_a-e_d)\tens
(qe_a-e_d)\over [2]_q}+q(q-1)e_a\tens e_a+\lambda
\theta\tens\theta\] is nondegenerate for $\lambda\ne
q(1-q)/[4]_q$, $\Delta_L$-invariant and symmetric in the sense
$\wedge(\eta)=0$
\end{lemma}
\proof This is adapted from \cite{Ma:rieq} and its properties then
verified directly in our case. \eproof

We next define the antisymmetrization tensor by
\[ \eps_{ijkl}{\rm Top}=e_i\wedge e_j\wedge e_k\wedge e_l\]
where ${\rm Top}=e_a\wedge e_b\wedge e_c\wedge e_d$ is
$\Delta_L$-invariant and a basis of $\Omega^4$. We can then define
define \begin{eqnarray*}
\star(e_i)&=&d_1^{-1}\eps_{ijkl}\eta^{jm}\eta^{kn}\eta^{lp} e_p
\wedge e_n\wedge e_m\\
\star(e_i\wedge
e_j)&=&d_2^{-1}\eps_{ijkl}\eta^{km}\eta^{ln}e_n\wedge e_m\\
\star(e_i\wedge e_j\wedge
e_k)&=&d_3^{-1}\eps_{ijkl}\eta^{lm}e_m\end{eqnarray*} for some
normalisations $d_i$ to be chosen.

Note that all constructions here are $\Delta_L$-covariant, under
which the space $\Lambda^1$ is a direct sum
\[ \Lambda^1={\rm sl}_{2,q}\oplus \C\theta,\quad {\rm sl}_{2,q}
\equiv\{e_b,e_c,e_z\equiv qe_a-q^{-1}e_d\}.\] and $\star^2$ has
these as eigenspaces. We now adjust $\lambda$ so that the
associated eigenvalues are the same.

\begin{lemma} For all invertible $q^2\ne 1$ there exists
precisely one value, $\lambda=q(1-q-q^2)/[2]_q$, such that $\eta$
is invertible and $\star^2\propto\id$ on $\Lambda^1$. In this case
we suppose $[3]_q\ne 0$ and set
\[ d_1=2q^2(1-q+q^2)[3]_q,\quad d_2=q^2[2]_{q^2},\quad d_3=q^2\]
Then, \begin{eqnarray*} &&\star e_a=-e_{abc}-\mu e_{bcd},\quad
\star e_b=-e_{abd},\quad
 \star e_c=q^2 e_{acd},\quad \star e_d=e_{bcd}\\
 && \star e_{ab}=-e_{ab}+2\mu e_{bd},\quad \star
 e_{ac}=e_{ac},\quad \star e_{ad}={1\over
 [2]_{q^2}}(2e_{bc}-q^2\mu e_{ad})\\
 &&\star e_{bc}={q^2\over [2]_{q^2}}(2e_{ad}+\mu e_{bc}),\quad
 \star e_{bd}=e_{bd},\quad \star e_{cd}=-e_{cd}\\
 &&\star e_{abc}=-e_a-\mu e_d,\quad \star e_{abd}=-e_b,\quad \star
 e_{acd}=q^{-2}e_c,\quad \star e_{bcd}=e_d
 \end{eqnarray*}
and $\star^2=\id$ on all degrees. The spaces of self-dual and
antiselfdual 2-forms are each 3-dimensional. We define $\star$
directly by these formulae for all invertible $q^2\ne \pm 1$.
\end{lemma} \proof We first compute $\eps$ as defined above. Its
nonzero values are
\begin{eqnarray*}
&&\eps_{1141}=-\eps_{1114}=-\eps_{1312}=\eps_{1411}=-\eps_{1414}
=\eps_{3121}=-\eps_{4111}=\eps_{4141}=
-q^2\eps_{1213}=q^2\eps_{2131}=\mu\\
&&\eps_{1234}=-\eps_{1243}=-\eps_{1324}=\eps_{1342}=\eps_{1423}
=-\eps_{1432}=-q^2\eps_{2134}=q^2\eps_{2143}=
\eps_{2314}=-\eps_{2341}\\
&&=-q^2\eps_{2413}=\eps_{2431}=q^{-2}\eps_{3124}=-q^{-2}\eps_{3142}
=-\eps_{3214}
=\eps_{3241}=q^{-2}\eps_{3412}=-\eps_{3421}=-\eps_{4123}\\
&&=\eps_{4132}=q^2\eps_{4213}=-\eps_{4231}=-q^{-2}\eps_{4312}
=\eps_{4321}=1. \end{eqnarray*} Using this, we define $\star$
(without normalisations) and compute $\star^2$ on $\Lambda^1$. We
solve for $\lambda$ such that its two eigenvalues coincide. This
has one solution which is such that $\eta$ is degenerate, and the
one shown. We then find that $\star^2\propto\id$ in degree 2 also,
and normalise $\star$ so that $\star^2=\id$ in all degrees. This
only fixes the product $d_1d_3$ but we chose these to reduce
repeated factors in $\star$. Also, it is clear by inspection that
\[ \Lambda^2_+=\{e_{bd},e_{ac},e_{ad}+e_{bc}\},
\quad \Lambda^2_-=\{e_{cd},e_{ab}-\mu
e_{bd},e_{ad}-q^{-2}e_{bc}\}.\] \eproof

Note that for the special value of $\lambda$ found in the
proposition above, which we use from now on, we have
\begin{eqnarray*}\eta&=&e_b\tens e_c+q^2 e_c\tens e_b-q^2(e_a\tens
e_d+e_d\tens
e_a+\mu e_d\tens e_d)\\
&&=e_b\tens e_c+q^2 e_c\tens e_b+q^2{(1-q^2)\over
[2]_{q^2}}e_z\tens e_z-{q^4\over[2]_{q^2}}\theta\tens\theta.
\end{eqnarray*}
This is precisely (in some conventions) the metric of
$q$-Minkowski space with $\theta$ the time direction. Likewise
$\eps$ is basically that for exact differentials on $q$-Minkowski
space in that context, see \cite{Ma:book}. In our case however,
the space is $SU_q(2)$ so there is no `time coordinate'. Instead,
$\theta$ being a generator of $H^1$, we see that the `time
direction' is created by q-deformation of the differential
calculus on $SU_2$ but is not exact, i.e. not $\extd$ of any time
coordinate.

With these general-$q$ preliminaries, we specialise from now on to
the reduced quantum group at the root given by $r=3$. We obtain
all specific formulae for this case, but expect similar features
for all odd $r$ as discussed at the end. We actually obtain such
results in the basis $\{a^mb^nc^k\}$ whereas the natural answers
equally involve the variable $d=a^2(1+q^2bc)$, to which we convert
using the identities
\begin{eqnarray*}
&&d^2=a(b^2c^2-qbc+1),\quad d^2b=-q(ab^2c-q^2ab),\quad d^2c
=-q(abc^2-q^2ab)\\
&&db^2=a^2b^2,\quad dc^2=a^2c^2,\quad
d(bc-q)=q^2(a^2b^2c^2-q^2a^2).
\end{eqnarray*}

We say that a form is {\em harmonic} if it is closed and coclosed.
The latter means in the kernel of
$\delta=\star\circ\extd\circ\star$. Likewise, coexact means with
respect to $\delta$, i.e. that the Hodge $\star$ of the form is
exact.

\begin{propos} At least for $r=3$, the element $\theta$ is coexact.
The element $\star\Theta$ is not closed. Moreover, $H^1$, $H^2$
have a basis of harmonic representatives, while the space of
harmonic elements of $H^3$ is the 3-dimensional kernel of
$\theta\wedge$.
\end{propos}
\proof  For $H^1$ the first three representatives are already
harmonic, while $n$ can be replaced by a harmonic 1-form
\[ h_3=q e_z-q^2 e_bd^2b +e_c a^2c.\]
One can also put $-qe_a$ for the first term since the difference
is $\theta$ already in the basis. For $H^2$ the $m_1,\cdots, m_4$
are already harmonic since they are up to a linear combination
self-dual or antiself-dual. They become part of our harmonic
(anti)self-dual basis
\[ h_1^+=e_{bd}ac^{r-1}, \quad h_1^-=(e_{ab}-\mu e_{bd})ac^{r-1},
\quad h_2^+=e_{ac}a^{r-1}b^{r-1}, \quad h_2^-
=e_{cd}a^{r-1}b^{r-1}.\] The remaining $m_5,m_6$ can be replaced
by harmonic ones
\[ h_3^+=e_{bd}d^2b+e_{ac}a^2c-(e_{ad}+e_{bc}),\quad
h_3^-=qe_{cd}a^2c+(e_{ab}-\mu e_{bd})d^2b+(e_{ad}-q^{-2}e_{bc})
\]
which are respectively self-dual and antiself-dual. The facts on
$\theta,\Theta$ can be directly verified. Finally, we take a basis
of Harmonic 3-forms and eliminate all those that are exact. This
leaves only three. Hence the dimension of the quotient is at most
3. On the other hand three harmonic 3-forms linearly independent
in the quotient are provided by applying $\star$ to the above
harmonic representatives of $H^1$. Up to coboundary and
normalisation, this gives a basis by $\Theta,h_1^*,h_2^*$ as
before and
\[h_3^*=e_{abd}d^2b+e_{acd}a^2c-(e_{abc}+e_{bcd}).\]
We can also write $qe_{abc}$ for the last term here since the
difference is a multiple of $\star\theta$ and this is exact.

Note also that in this harmonic basis the action of $\theta\wedge$
in Proposition~3.5 is more symmetric. It clearly sends harmonic
forms to harmonic forms. In fact we find $\theta\wedge
h_1=h_1^--q^{-2}h_1^+$, $\theta\wedge h_2=h_2^+-h_2^-$ as before
and $\theta\wedge h_3=h_3^+-q^2h_3^-$. Their complement has basis
$h_1^++q^2h_1^-$, $h_2^++h_2^-$ as before, and $h_3^++q^2h_3^-$.
The image of these under $\theta\wedge$ is now a multiple of
$h_1^*,h_2^*,h_3^*$, with complement $\Theta$. \eproof

This immediately implies that there is no Hodge decomposition
theorem (into a direct sum of exact, coexact and harmonic forms in
each degree), precisely because $\theta$ is a nonzero element that
is both coexact and harmonic.

Leaving now cohomology, we consider general forms and `wave
equations'. As well as the operator $\extd+\delta$ who's kernel is
the harmonic forms (given that they map into different degrees),
we also have the Laplacian $\square=\delta\extd+\extd\delta$.

\begin{propos} For $r=3$ the dimensions over $\C$ of the spaces of
Harmonic forms and the kernel of $\square$ are shown in Table 1.
Also for comparison we remind the dimensions of the closed and
exact forms in each degree as found in Section~3. Coclosed and
coexact are given by reversing the relevant rows. In particular,
${\rm harmonic}\subset\ker\square$ is strict.
\end{propos}
\proof This is direct computation once the matrices for the
various operators above have been found explicitly. \eproof
\begin{table}
\[\begin{array}{c|ccccc}
r=3& \Omega^0& \Omega^1& \Omega^2 &\Omega^3 &\Omega^4 \\
\hline {\rm All}& 27& 108& 162& 108& 27\\
{\rm Closed} & 1& 30 &84 & 82&27 \\
{\rm Exact} & 0 & 26 & 78 & 78& 26 \\
{\rm Harmonic} &1 &16&30&16&1   \\
\ker\square &13& 33 & 40 & 33 & 13  \\
\end{array}\]
\caption{Number of independent forms of various types in each
degree, for $r=3$.}
\end{table}

Next we look in detail at physical `wave equations'. For spin 0 or
scaler fields, we find that $\square$ is {\em not} fully
diagonalisable. This is related to the nonsemisimplicity of the
Hopf algebra.

\begin{propos} For $r=3$ a full set of 13 zero-modes of $\square$
in spin zero are
\[1,\ a,\ b,\ c,\ d,\ ab^2,\ a^2b,\ db^2,\ d^2b,\ ac^2,\ a^2c,\
dc^2,\ d^2c.\]
  In
addition there are 9 `massive' modes of eigenvalue $6(q+1)$ given
by
\[ a^2,\ b^2,\ c^2,\ d^2,\ ab,\ ac,\ db,\ dc,\ bc-1.\]
\end{propos}
\proof Elementary computation once $\square$ is defined. Note the
zero modes $ab^2c-q^2ab$ and $a^2c$ already featuring in the
construction of harmonic forms above. \eproof

Note that we do not consider `orthogonality' since the correct
`reality' properties are not clear when $q$ is a root of unity.
Instead we are guided at our algebraic level by simplicity of
expressions. It is worth noting that there is, however,
necessarily a translation-invariant `integral' functional in the
Hopf algebra sense.

Next we solve the `spin 1' or 1-form system. Following the
notations in physics, we say that a 1-form is in {\em Lorentz
gauge} if it is coclosed. It is in {\em temporal gauge} if it can
be written entirely in terms of $e_b,e_c,e_z$ (i.e. no $\theta$
component when taking these four as basis). By number of `modes'
we will mean only the dimensions of the relevant spaces or
quotient spaces (the number of linearly independent vectors in any
basis).

\begin{propos} Let ${\rm Max}=\delta\extd$ be the Maxwell operator
on $\Omega^1$. Then for $r=3$: (i)  $\ker{\rm Max}$ is
54-dimensional, hence up to gauge equivalence (i.e. modulo exact
1-forms) there are 28 `true' spin 1 zero modes, of which exactly 4
have zero curvature $\extd A$ (namely the harmonic basis of
$H^1$). (ii) If we `gauge fix' to Lorentz gauge by looking among
coclosed 1-forms $A$ then there are 32 zero modes but only 20 true
ones when taken modulo exact. (iii) Ditto for temporal gauge. (iv)
Every zero mode is gauge equivalent to the sum of a zero mode in
Lorentz gauge and one in temporal gauge (with 12 modes in both
gauges up to equivalence.)
\end{propos}
\proof We compute the dimension of the kernel of ${\rm Max}$ as
54. It contains the exact 1-forms, so subtracting 26 gives the
true dimension `modulo gauge'.  Much more work gives explicit
bases of representatives of the various types of modes constructed
as kernels of suitable linear maps. Here we use the same method as
for the cohomology computations, namely we first eliminate all
elements of the relevant kernel which are exact. The remainder
could still be linearly dependent in the quotient. We then
painstakingly chose enough representatives to give the required
dimensions, i.e. checking that together with the image of
$\extd_0$ they form a basis of the original kernel. Of course this
process is not unique (we choose the simplest representatives
where possible). Explicitly, they are as follows. 20 modes obeying
the gauge fixing are the elements $\{h_1,h_2,h_3\}$ of the
Harmonic basis of $H^1$ above (with zero curvature) plus the 13
modes of the form
\[ A=\theta f,\qquad \square f=0\]
(i.e. induced by spin 0 solutions as given above), and 4 more
coclosed modes which we have to specify. E.g. the vector space of
coclosed modes which are also in temporal gauge is 19 dimensional,
reducing to 7 true modes in the quotient, of which 3 are
$\{h_1,h_2,h_3\}$ already counted. The remaining four are:
\begin{eqnarray*} A_1&=&e_z d(bc-q)- e_bb(bc-q^2)+q e_cd^2c\\
A_2&=& qe_zabc -e_ba^2b+e_cdac \\
A_3&=&qe_zb^2c-qe_bab^2+e_c d^2a
  \\
A_4&=&e_zbc^2 -qe_ba^2d+e_cdc^2.
\end{eqnarray*}
Finally, we must complete the basis with 8 modes which are not,
however, coclosed. We find that the dimension of the space of
solutions in temporal gauge is (like for Lorentz gauge) 32
dimensional, reducing to 20 true temporal-gauge modes in the
quotient. We have already used 7 of them above and can choose 8
more from among the remainder, e.g.
\begin{eqnarray*}
&& A_5=e_b a^2c,\quad A_{6}=e_cd^2b,\quad A_7=qe_c
b^2c-e_zab^2,\quad
A_8=qe_b+e_za^2c,\quad A_9=q^2e_bab^2+e_zdc^2\\
&& A_{10}=e_c abc-qe_z a^2b,\quad A_{11}=e_bdb^2+q^2e_zd^2b,\quad
A_{12}=e_ba^2-qe_zd^2c.
\end{eqnarray*}
We can also use $A_6'= e_zdb^2$ in place of $A_6$ with the same
curvature up to normalisation. In this way we may `patch' the
moduli of solutions into Lorentz and temporal gauge, with some
overlap. We could equally chose 20 temporal gauge modes and
complete with 8 more in Lorentz but not temporal gauge if we
preferred. This means that there are 12 true modes which can be
viewed in either gauge by a gauge transformation (but only 7 which
can be transformed to a solution in both gauges simultaneously as
explained above). \eproof

These results show several key features of the electromagnetic
theory. First and foremost, there are 24 `electromagetic' modes
with nonzero curvature $F=\extd A$ obeying the source-less Maxwell
equation (i.e. forming a basis with the zero curvature ones). They
are the analogue of the photon self-propagation modes in usual
physics. I.e. `there is light'. The remaining 4 modes of zero
curvature indicate nontrivial topology and the existence of the
`Bohm-Aharanov' effect. Finally, we see that usual gauge fixing to
`Lorentz gauge' (where $\delta A=0$) does not work: not all
solutions obey the gauge fixing condition. Likewise for temporal
gauge fixing. Such problems can potentially plague any nontrivial
gauge theory but here in our concrete model we see how the moduli
space can instead be `covered by patches' built from Lorentz and
temporal gauge. Note also that two representatives in Lorentz
gauge can only differ by $\extd f$ with $f\in\ker\square$ and in
usual electromagnetism this would be forced to be zero by boundary
conditions at infinity (so that there would be a unique
representative fixed by the gauge condition); in our case we do
not have any such natural conditions, i.e. the possibility of
nontrivial `Gribov ambiguities'\cite{ItzZub}. This would be
relevant to the quantum electromagnetic theory if one tried to
impose gauge fixing in the functional integral.

\begin{propos} Of the 28 true zero modes of $\rm Max$ for $r=3$,
exactly 16 have self-dual curvature and 16 antiself-dual
curvature. Every zero mode is gauge equivalent to the sum of a
self-dual and an antiself-dual zero modes (with the four zero
curvature modes in both classes).
\end{propos}
\proof We compute dimensions as kernels of suitable maps. Thus the
space of 1-forms with self-dual curvature is 42 dimensional,
reducing to 16 true modes in the quotient. Similarly for
antiself-dual. We next proceed to find reasonable representatives
in the space of 1-forms modulo exact ones forming a basis as in
the computations above. To this end, we also note that the
dimension of the space of 1-forms which have self-dual curvature
and are coclosed is 20-dimensional, reducing to 8 true modes in
the quotient, including the 4 of zero curvature (the harmonic
basis of $H^1$) already given. This leads us to 8 of our basis of
16 forms given by $\{\theta,h_1,h_2,h_3\}$ and
\[A_1=e_aa,\quad A_2=e_ab,\quad A_3=e_d c,\quad A_4=e_dd\]
say (these are equivalent up to normalisation and coboundaries to
coclosed modes but not themselves coclosed). We complete the basis
of 1-forms with self-dual curvature by
\begin{eqnarray*}
&& A_5=(\mu e_d+e_a)a^2b+e_c abc,\quad  A_6=(\mu
e_d+e_a)ab^2+q^2e_cb^2c\\
&& A_7=e_da^2c^2+q^2e_b bc^2,\quad A_8=e_dd^2c+q e_bdbc\\
&& A_9=e_adb^2,\quad A_{10}=e_d ac^2,\quad
A_{11}=e_b-e_aa^2c,\quad A_{12}=e_ad^2b-qe_b dc^2.
\end{eqnarray*}
These are all chosen to have simple expressions for their
self-dual curvatures, namely (with $e_+\equiv e_{ad}+e_{bc}$, and
up to normalisations),
\begin{eqnarray*}
&&F_1=e_+a+q^2 e_{ac}c,\quad  F_2=e_+b+q^2 e_{ac}d,\quad
F_4=e_+d+q e_{bd}b,\quad F_3=e_+c+q e_{bd}a\\ && F_5=e_{ac}a,\quad
F_6=e_{ac}b,\quad F_7=e_{bd}c,\quad F_8=e_{bd}d\\
&& F_9=e_{ac}d^2b-q e_+db^2,\quad
F_{10}=e_+ac^2-qe_{bd}a^2c\\
&& F_{11}=e_+a^2c-q^2 e_{ac}ac^2-e_{bd},\quad
F_{12}=e_{ac}+qe_{bd}db^2-e_+d^2b.
\end{eqnarray*}
These are exact and coclosed (hence harmonic) 2-forms. Note that
this is possible because the Hodge decomposition again does not
hold, here in degree 2. We can similarly find 12 antiself-dual
forms completing the zero curvature ones to a basis of the
antiself-dual moduli space. One then checks that these 12, the
above 12 self-dual modes and the 4 zero curvature modes are
linearly independent modulo exact forms.  This decomposition also
holds before working modulo exact forms, with the 30 closed forms
as the intersection of the two 42-dimensional spaces. \eproof

It turns out that we can also `patch' the moduli of solutions of
the sourceless Maxwell equations into Lorentz gauge and self-dual
ones. Here the self-dual modes $A_5,\dots, A_{12}$ are beyond the
reach of the Lorentz gauge fixing condition, being linearly
independent modulo exact forms to the basis of the Lorentz
gauge-fixed solutions in Proposition~4.6. Similarly for temporal
gauge.

\begin{corol} At least for $r=3$, (i) every zero mode of ${\rm Max}$
is gauge equivalent to the sum of one of the form $\theta f$ where
$\square f=0$ and a self-dual one (with the mode $\theta$ in both
classes). (ii) every zero mode of $\rm max$ is gauge equivalent to
the sum of one in temporal gauge and a self-dual one (with 8 modes
including 4 zero curvature ones in the overlap).Similarly using
antiself-dual modes.
\end{corol}\proof (i) We check that the 16 modes $\theta f$,
$\square f=0$ and $\{h_1,h_2,h_3\}$ in
the proof of Proposition~4.6 are linearly independent modulo exact
forms from the $A_1,\cdots, A_{12}$ self-dual modes in
Proposition~4.7. Also note that if we want to have as much as
possible of the basis in Lorentz gauge then we could equally well
use the coclosed self-dual modes
\begin{eqnarray*} A'_1&=&e_a a-e_ba^2b+q e_c(bc-q)c+qe_zabc \\
A'_2&=&q e_ab-\theta b+qe_zb^2c-e_bab^2+e_cd(bc-q) \\
A'_3&=&
e_d c+qe_zbc^2-qe_babc+qe_ca^2c^2 \\
A'_4&=&e_d d-e_z(d(bc-q)-q^2c)-qe_cd^2c+e_b(b^2c+qa)
.\end{eqnarray*} These are gauge equivalent (up to normalisation)
to the $A_1,\dots, A_4$ in Proposition~4.7, giving a full set of
20 coclosed modes and a basis along with the $A_5,\cdots,A_{12}.$
(ii) For temporal gauge we find that there is similarly a
20-dimensional space of forms which are both in temporal gauge and
have self-dual curvature, reducing to 8 in the quotient. They
include 4 of zero curvature (so $H^1$ has a basis of
representatives in temporal gauge) and all 8 are in fact gauge
equivalent to the above 8 modes that were self-dual and
renderabale in Lorentz gauge (the $\{\theta,h_1,h_2,h_3\}$ and
$A_1,\cdots, A_4$ (or $A'_1,\cdots, A'_4$) as just discussed). So
the self-dual $A_5,\cdots,A_{12}$ in Proposition~4.7 again
complete to a full set of self-dual forms. This time we find that
the 12 temporal gauge modes $A_1,\cdots A_{12}$ in Proposition~4.6
then complete to full set of 28 zero modes of $\rm Max$. \eproof

Finally, we give examples of some source $J$ and solve the full
Maxwell equation  $\delta F=J$. Recall that the element $\theta$
is coexact as one would need for any source $J$.

\begin{propos} For $r=3$, a basis of valid sources (i.e. in the
image of $\rm Max$) in the direction of $\theta$ is provided by
\[ \theta\{1,a,b,c,d\}.\]
In particular, the element $\theta$ is a valid source and has a
gauge field (not uniquely determined since we can add any of the
above zero modes) given in Lorentz gauge by
\[ A=-{q^2\over 12}\theta\, bc(1+bc)-{q\mu\over 12}(e_a+e_c a^2c).\]
Its curvature is
\[ F={q\over 4}e_{ad}-{\mu\over 12}\left((e_{ab}-e_{bd})d^2b
+q(e_{cd}-e_{ac})a^2c\right)\]
\end{propos}
\proof We first compute the dimension of the subspace of the image
of $\rm Max$ of the form $\theta f$ as 5. This is found as the
dimension of the image of $\rm Max$ minus that of the image of
$T\circ {\rm Max}$ where $T$ is the linear map whose kernel is
spanned by $\theta$ over $\CA$. We then solve explicitly for the
example $J=\theta$ in Lorentz gauge. Note that the second term in
$A$ is topological, being a multiple of the fourth basis element
of $H^1$ in Section~3. It can be omitted so that $A$ is itself
$\theta$ times a function (if we abandon the Lorentz gauge),
without changing the curvature. \eproof

According to the physical picture mentioned above, $\theta$ could
be viewed as a Minkowski time direction. So there is a `current'
in the cotangent space of $SU_q(2)$ in this direction (but no
actual current flow as time is not a coordinate) generating this
gauge field. In usual Maxwell theory such a current in the time
direction corresponds to a static electric charge density.
Accordingly, the source $\theta$ can be viewed as a uniform charge
density over the noncommutative $S^3$ leading to gauge field and
(electric) curvature field as stated. There are of course many
other sources, the dimension of the image of $\rm Max$ being 54
(for $r=3$).

\begin{propos} For $r=3$, the subspace of sources in the `spatial'
directions spanned by $e_z,e_b,e_c$ is 40-dimensional. Those
purely along each of the three directions have bases
\[ e_z,\quad e_b\{1,c^2,d^2,dc,dc^2,d^2c\},\quad
e_c\{1,a^2,b^2,ab,ab^2,a^2b\}.\]
In particular, the gauge fields and their curvatures
\begin{eqnarray*}&A_1={q^2\over 6}e_z,\quad\quad &F_1={q^2\mu\over
6}e_{bc}\\
&A_2={q^2\over 6}e_b,\quad\quad &F_2=-{\mu\over 6}( q^2
e_{ab}+e_{bd})\\
 &A_3=-{1\over 6}e_bdb^2,\quad \quad &F_3=-{\mu\over 6}\left(
 e_{bc}d^2b+(q^2 e_{ab}+e_{bd})db^2\right)
\end{eqnarray*}
are solutions for the sources $e_z,e_b,e_c$ respectively.
\end{propos}
\proof Here we compute the dimension of the subspace of spatial
currents as the image of $\rm Max$ as the dimension of the image
of $\rm Max$ minus that of the image of $S\circ {\rm Max}$ where
$S$ is the linear map whose kernel is the spatial directions.
Similarly along each of the directions $e_z,e_b,e_c$ separately we
obtain dimensions 1,6,6 respectively. We then find the right
number of independent modes. Finally, we solve for some of these
in temporal gauge and exhibit the solutions for the three constant
sources directions. Again, these solutions are not unique since we
can add any of the above zero modes of $\rm Max$. \eproof

In classical electrodynamics spatial sources would correspond to
currents inducing magnetic configurations. Another example is the
source $e_c b^2$ having a solution with curvature proportional to
$(e_{cd}+qe_{ac})b^2$. From these various (and other) solutions we
see that the natural electric and magnetic curvature directions
under this time/space decomposition are spanned by
\[ \Lambda^2_E=\{e_{ad},\, e_{ab}-e_{bd},\, e_{cd}-e_{ac}\},\quad
\Lambda^2_B=\{e_{bc},\, q^2e_{ab}+e_{bd},\,e_{cd} +q e_{ac}\}\]
respectively.

\begin{table}
\[\begin{array}{c|ccccc}
r=5& \Omega^0& \Omega^1& \Omega^2 &\Omega^3 &\Omega^4 \\
\hline {\rm All}& 125 & 500 & 750 & 500 & 125 \\
{\rm Closed} & 1& 128 & 378 & 376& 125 \\
{\rm Exact} & 0 & 124 & 372 & 372 & 124 \\
{\rm Harmonic} &1 &36&70&36&1   \\
\end{array}\]
\caption{Number of independent forms of various types in each
degree, for $r=5$.}\end{table}
\begin{table}
\[\begin{array}{l|cccc}
{\rm Maxwell}& r=3 & & r=5  \\
\hline {\rm All\ zero\ modes}& 28 & (54)& 68 &(192) \\
{\rm Coclosed} & 20 & (32)& 52 & (84)   \\
{\rm Temporal} & 20 & (32) & 52 & (84)  \\
{\rm Cocl.}\cap{\rm Temp.} & 7 & (19) & 19& (51)  \\
{\rm self-dual}& 16 &(42) & 36 &(160) \\
{\rm zero\ curv.}& 4 & (30)& 4 & (128)\\
{\rm Cocl.}\cap{\rm s.d.}& 8 & (20)& 20 & (52)\\
{\rm Temp.}\cap{\rm s.d.}& 8 & (20)& 20 & (52)\\
\theta f\ {\rm modes}& 13 &(13) &33& (33)\\
 \hline
{\rm All\ sources}&  54 & & 308\\
{\rm spatial\ sources}& 40 & & 216\\
\theta f\ {\rm sources} & 5& & 17\\
\end{array}\]
\caption{Summary of electromagnetic theory for $r=3,5$. Number of
independent solutions of the sourceless Maxwell equations modulo
exact forms (in brackets before making the quotient). We also
summarize the types of valid sources.}
\end{table}

It should be mentioned in conclusion that other odd roots appear
to give similar features as the $r=3$ case. The preliminary
Table~2 summarizes the form dimensions for $r=5$, after which
Table~3 summarizes the Maxwell theory above and the corresponding
numbers for $r=5$. From these and further inspection we find the
same qualitative features, e.g. all solutions can be written as
sums of self-dual and antiself-dual solutions with overlap given
by the zero curvature modes (both modulo exact 1-forms and before
taking the quotient); the temporal and Lorentz gauges patch the
moduli space, etc. Also, $\theta,e_z,e_b,e_c$ are all valid
sources of electric and magnetic type among others in the numbers
shown. From the tables we note another novel feature that we
expect for all odd $r$, namely a linear isomorphism between
harmonic 1-forms and self-dual solutions modulo exact forms.

\bigskip

\baselineskip 14pt

\end{document}